\documentclass{amsart}
\usepackage[utf8]{inputenc}

\title{A note on the Nielsen realization problem for Hyper-Kähler manifolds}
\author[S.\ Billi]{Simone Billi}
\address{University of Genova, Department of Mathematics, Genova 16146, Italy}
 \email{{\tt simone.billi@edu.unige.it}}

\usepackage[english]{babel}
\usepackage[T1]{fontenc}
\usepackage{amssymb,amsmath,amsthm,amsfonts}
\usepackage{thmtools,mathtools}
\usepackage{microtype}
\allowdisplaybreaks
\usepackage{enumitem}
\usepackage{tikz,tikz-cd}
\usepackage[breaklinks,pdfencoding=auto,psdextra]{hyperref}
\usepackage{bm} 
\usepackage[scaled]{beramono}

\usepackage{listings}
\usepackage{caption}
\usepackage{xcolor}
\usepackage{amsmath}
\usepackage{amssymb}
\usepackage{amsfonts}
\usepackage{mathrsfs}
\usepackage{amsthm}
\usepackage{graphicx}
\usepackage{floatflt}
\usepackage{epstopdf}
\usepackage{comment}
\usepackage{tikz-cd}
\usepackage{url}
\usepackage{appendix}

\lstdefinelanguage{Macaulay2}{
comment=[l]{--},
alsoletter={'},
alsoother={_},
}
\DeclareMathOperator{\bL}{\mathbf{L}}
\DeclareMathOperator{\bLambda}{\mathbf{\Lambda}}

\DeclareMathOperator{\Aut}{Aut} 
\DeclareMathOperator{\og10}{OG10}
\DeclareMathOperator{\Isom}{Isom} 
 
\DeclareMathOperator{\Homology}{H}
\DeclareMathOperator{\Gr}{Gr}

\DeclareMathOperator{\T}{T}

\DeclareMathOperator{\Diff}{Diff}
\DeclareMathOperator{\Mod}{Mod}
\DeclareMathOperator{\Mon}{Mon}

\DeclareMathOperator{\E}{E}
\DeclareMathOperator{\B}{B}
\DeclareMathOperator{\M}{M}

\DeclareMathOperator{\id}{id}
\DeclareMathOperator{\SO}{SO}

\DeclareMathOperator{\bO}{\mathrm{O}}
\DeclareMathOperator{\II}{\mathbf{I}}

\newtheorem{thm}{Theorem}[section]
\newtheorem{lem}[thm]{Lemma}

\newtheorem{definition}[thm]{Definition}
\newtheorem{cor}[thm]{Corollary}

\newtheorem*{problem*}{Problem}
\newtheorem{question}{Question}
\newcommand{\hk}{hyper-K{\"a}hler}

\newcommand{\KTST}{K3$^{[n]}$ type}

\newtheorem{rmk}{Remark}
\begin{document}
\maketitle

\begin{abstract}
    We give an answer to the Nielsen realization problem for \hk{} manifolds in analogy to the case of K3 surfaces. We determine that, for some of the known deformation types, the representation of the mapping class group on the second cohomology admits a section on its image. 
\end{abstract}
\section*{Introduction}
The Nielsen realization problem was originally formulated by Nielsen in \cite[Section 4]{nielsen1942abbildungsklassen}, and then affirmatively solved in \cite{kerckhoff1983nielsen}. The question is whether any finite group $G$ of mapping classes of a complex curve can be lifted to a group of diffeomorphisms (which preserve the metric and the complex structure). An equivalent formulation of the problem is if $G$ fixes any point in the Teichm\"uller space.  \\
The answer to this question for K3 surfaces is given by Farb and Looijenga in \cite{farb2021nielsen}: their result shows that, if \(S\) is a K3 surface, then not every finite subgroup of the mapping class group $G\subset\Mod(S)=\pi_0(\Diff^+(S))$ can be lifted, but there is a \(G\)-invariant lattice $\bL_G$ which determines if it is possible or not. 
The lattice \(\bL_G\) is a sublattice of the second cohomology lattice of the K3 surface. In relation to this,
Baraglia and Konno proved in \cite{Nielsen} that for a K3 surface the representation of the mapping class group on the second cohomology group has a section over its image.
\\

We address the same problems for \hk{} manifolds. We consider the lattice \(\bLambda_X=(\Homology^2(X,\mathbb{Z}),q_X)\) where \(q_X\) is the Beauville-Bogomolov-Fujiki quadratic form, its group of orientation-preserving isometries \(\bO^+(\bLambda_X)\) and the group \(\Aut(X)\) of biholomorphic automorphisms of \(X\) with respect to a chosen complex structure. In the first section we prove the following:
\begin{thm}\label{thm1}
	Let \(X\) be an \hk{} manifold such that \(\Aut(X)\rightarrow \bO^+(\bLambda_X)\) is injective, then the restriction of  \(\rho\colon\Mod(X)\rightarrow \bO^+(\bLambda_X)\) to the stabilizer \(\Mod(X)_{\mathcal{C}}\) of a component \(\mathcal{C}\) of the Teichmüller space has a section over its image \(\Mon^2(X)\).
\end{thm} The group \(\Mon^2(X)\) has finite index in \(\bO^+(\bLambda_X)\) by \cite[Theorem 2.6 (i)]{verbitsky2020errata}. The previous result applies to certain classes of \hk{} manifolds:
\begin{cor}\label{cor:max_monodromy}
	If \(X\) is \KTST{} with \(n-1\) a power of a prime, or \(X\) is of \(\og10\)-type, then \(\rho\colon\Mod(X)\rightarrow \bO^+(\bLambda_X)\) has a section.
\end{cor}

In the second section we define the lattice $\bL_G$ for a finite subgroup \(G\subset \Mod(X)\) and \(X\) a \hk{} manifold in analogy to the case of K3 surfaces, see \autoref{def:LG} and the previous lines. Let \(\mathcal{T}_{Ein}\) be the Teichmüller space of Einstein metrics on \(X\) and \(\Delta_{\mathcal{C}}\) the set of negative forms on \(\bLambda_{X}\) associated to rational curves for a metric in the connected component \(\mathcal{C}\). We prove the following:
\begin{thm}\label{thm3}
	Let \(X\) be a \hk{} manifold and let $G$ be a finite subgroup of $\Mod(X)$.
	\begin{itemize}
		\item $G$ lifts to a group of isometries of an Einstein metric if and only if $G$ fixes a connected component $\mathcal{C}$ of $\mathcal{T}_{Ein}$ and $\bL_G$ does not contain any element of $\Delta_\mathcal{C}$.
		\item $G$ lifts to a group of isometries for a \hk{} metric if and only if $G$ fixes a connected component $\mathcal{C}$ of $\mathcal{T}_{Ein}$, $\bL_G$ does not contain any element of $\Delta_\mathcal{C}$ and $\bL_G^\perp$ contains the trivial representation (in this case the metric can be chosen so that $X$ is projective and $G$ acts by algebraic automorphisms). 
	\end{itemize}
	Similarly, a finite subgroup of $\bO^+(\bLambda_X)$ which is contained in $\Mon^2(X)$ is the image via \(\rho\) of a finite subgroup of \(\Mod(X)\) that lifts under the same conditions.
\end{thm}

\section*{Acknowledgments}
The author would like to thank Giovanni Mongardi for suggesting to tackle these problems, for his useful comments and the supervision of the project. Moreover, the author thanks Kejia Zhu for pointing out a gap in the previous version of this note, Nikola Sadovek for discussing \autoref{homotopy_equiv} and the referee for useful comments about \autoref{thm1} and its proof.

The author was partially supported by PRIN 2020KKWT53 003 - Progetto: Curves, Ricci flat Varieties and their Interactions. The author is member of the INdAM group GNSAGA, and was partially supported by it.

\section*{Setting}

We will denote by $X$ the underlying smooth manifold of a \hk{} manifold, i.e. a simply connected complex compact K\"ahler manifold such that $\Homology^0(X,\Omega^2_X)=\mathbb{C}\sigma$ where $\sigma$ is a non-degenerate holomorphic $2$-form. A K\"ahler class on $X$ determines a unique K\"ahler-Einstein metric \(g\) and a $2$-sphere of complex structures for which the metric is still K\"ahler. Let $\Mod(X)=\pi_0(\Diff^+(X)))$ be its mapping class group and $\bLambda_X=(\Homology^2(X,\mathbb{Z}),q_X)$ the second cohomology group endowed with the Beauville-Bogomolov-Fujiki quadratic form $q_X$, which by the Fujiki relations is a topological invariant. We are going to use the notation $S$ when dealing with a K3 surface.\\
There is a map from $\Mod(X)$ to $O(\bLambda_X)$ whose image \(\Gamma\) sits in the index $2$ subgroup $\bO^+(\bLambda_X)$, consisting of orientation-preserving isometries. The subgroup \(\Gamma\) has finite index in $O(\bLambda_X)$ by \cite[Theorem 1.1]{verbitsky2020errata}. \\
Consider the following maps
\[\Diff^+(X)\to \Mod(X)\xrightarrow{\rho} \bO^+(\bLambda_X),\]
if \(S\) is a K3 surface then $\rho\colon\Mod(S)\rightarrow \bO^+(\bLambda_S)$ has a section by \cite[Theorem 1.1]{Nielsen}.
A partial generalization of this result for \(X\) \hk{} is given by \autoref{thm1}.

Let $G$ be a finite subgroup of $\Mod(X)$, one can ask the following:
\begin{problem*}[Nielsen realization]
	Does there exists an Einstein metric $g$ on $X$ such that $G$ is realizable as a subgroup of $\Isom(X,g)$? Can the metric $g$ chosen to be also K\"ahler?  
\end{problem*}
An answer to this question for K3 surfaces is given in \cite[Theorem 1.2]{farb2021nielsen}, in \autoref{thm3} we give a similar answer for \hk{} manifolds.
\\
As in \cite{looijenga2021teichmuller}, we consider the following Teichm\"uller spaces associated to \(X\) with the respective period maps (see \cite[({\(\dagger\)) at page 270}]{looijenga2021teichmuller}):
\[
\begin{tikzcd}
\mathcal{T} \arrow[d, "\mathcal{P}"'] & \mathcal{T}_{HK} \arrow[r] \arrow[d, "\mathcal{P}_{HK}"'] \arrow[l] & \mathcal{T}_{Ein} \arrow[d, "\mathcal{P}_{Ein}"] \\
\mathcal{D}                          & \mathcal{D}_{HK} \arrow[r] \arrow[l]                        & {\Gr^+(3,\bLambda_X\otimes \mathbb{R})}              
\end{tikzcd}
\]
where \(\mathcal{T}\) is the Teichm\"uller space of complex structures, $\mathcal{T}_{HK}$ is the Teichm\"uller space of \hk{} metrics with unitary volume and $\mathcal{T}_{Ein}$ is the Teichm\"uller space of Einstein metrics with unitary volume, the sets of such up to isotopy. The space $\mathcal{D}=\{[v]\in\mathbb{P}(\bLambda_X\otimes\mathbb{C})| q_X(v)=0, q_X(v,\overline{v})>0\}$ is the period domain of $\mathcal{T}$,  the space \(\mathcal{D}_{HK}\) consists of pairs \((z,r)\) where \(z\in\mathcal{D}\) and \(r\) is a ray in \(\bLambda^{1,1}_X\otimes\mathbb{R}\) and the Grasmannian of positive-definite \(3\)-planes $\Gr^+(3,\bLambda_X\otimes \mathbb{R}) $ is the period domain of $\mathcal{T}_{Ein}$. The vertical maps are the period maps associating to a metric the class $[\sigma_X]$ in the first case, the additional ray determined by a Kähler metric in the second case and the positively-oriented $3$-plane $P$ determined by the real and imaginary parts of $\sigma_X$ in the third case. We refer to \cite[Sections 3.1, 3.2]{looijenga2021teichmuller} for a discussion about the period domain and the period map of \(\mathcal{T}\), to \cite[Section 3.4]{looijenga2021teichmuller} for the period domains and the period maps of \(\mathcal{T}_{HK}\) and \(\mathcal{T}_{Ein}\). The period maps give to the Teichm\"uller spaces the structure of respectively complex and differentiable spaces, see \cite[Section 3.3]{looijenga2021teichmuller} for a discussion about \(\mathcal{T}\) and \cite[Section 3.4]{looijenga2021teichmuller} for a discussion about \(\mathcal{T}_{HK}\) and \(\mathcal{T}_{HK}\) (which is there denoted by \(\mathcal{T}_{\mathbb{H}}\)). By the observations before and after \cite[Corollary 3.8]{looijenga2021teichmuller}, the map \(\mathcal{T}_{HK}\rightarrow \mathcal{T}\) has contractible fibers and  $\mathcal{T}_{HK}\rightarrow \mathcal{T}_{Ein}$ simply forgets the choice of the complex structure so it is a locally trivial $2$-sphere bundle. As a consequence, there is a bijection between the sets of connected components of the three spaces. Moreover, by \cite[Corollary 5.4]{looijenga2021teichmuller} each connected component $\mathcal{C}$ of $\mathcal{T}_{Ein}$ is mapped diffeomorphically onto a simply connected subset of $\Gr^+(3,\bLambda_X\otimes \mathbb{R})$, the Grassmannian of oriented positive definite $3$-planes. In particular, \(\mathcal{C}\) is simply connected. \\
Denote by $\Mod(X)_\mathcal{C}$ the $\Mod(X)$-stabilizer of the connected component $\mathcal{C}$ and consider the \textit{Torelli group} $\T(X)=\ker(\rho)$.
The \textit{monodromy group} $\Mon^2(X)\subseteq \Gamma\subseteq \bO^+(\bLambda_X)$ is the image of $\Mod(X)_\mathcal{C}$ and can be identified with the group of isometries coming from parallel transport operators \cite[Theorem 7.2]{Verbitsky_2013}.\\

In the case of K3 surfaces, $\Mod(S)_\mathcal{C}$ maps isomorphically onto $\bO^+(\bLambda_S)=\Mon^2(S)$ via $\rho$ giving the isomorphism 

\begin{align*}
\Mod(S)\cong \T(S)\rtimes \bO^+(\bLambda_S)
\end{align*}
which implies Theorem \ref{thm1}.
In this particular case, $\T(S)$ permutes transitively the connected components of $\mathcal{T}$ implying that the moduli space of marked K3 surfaces $\mathcal{M}_{\bLambda_S}=\mathcal{T}/\T(S)$ is connected.\\
\begin{rmk}
	If $X$ is \hk{} of dimension bigger than $2$ then $\Mod(X)$ could be just an extension of $\T(X)$ and $\Gamma$, similarly $\Mod(X)_\mathcal{C}$ could be an extension of $\Mon^2(X)$ and $\T(X)\cap \Mod(X)_\mathcal{C} $, but by \cite[Remark 2.5]{verbitsky2020errata} the intersection $\T(X)\cap \Mod(X)_\mathcal{C} $ is always finite. Moreover, $\T(X)$ acts on $\pi_0(\mathcal{T})$ with finitely many orbits, each connected component has finite stabilizer and an element of $\T(X)$ which fixes an element $g\in \mathcal{T}$ fixes the entire connected component of $g$ (\cite[Theorem 3.1]{verbitsky2020errata}). In general, the moduli space of marked \hk{}  manifolds $\mathcal{M}_{\bLambda_X}=\mathcal{T}/\T(X)$ could have more connected components, but each one is simply connected. 
\end{rmk}
Rephrasing, we have $\Mon^2(X)\cong\Mod(X)_\mathcal{C}$ precisely when $\T(X)\cap \Mod(X)_\mathcal{C}$ is trivial and $\Mod(X)\cong \T(X)\rtimes \Gamma$ exactly when $\rho$ admits a section on its image. There could be a section of $ \rho$ over its image even if $\Mon^2(X)$ is a proper subgroup of $\bO^+(\bLambda_X)$ and on the other hand a priori there is still the possibility that $\Mon^2(X)=\Gamma=\bO^+(\bLambda_X)$ but $\Mod(X)$ is a not the semidirect product of $\T(X)$ and the stabilizer of a component.\\

\section{Sections of the representation map}
Here we generalize the proof of \cite[Theorem 1.1]{Nielsen} in some cases. Note that in this setting most of the groups we consider are discrete. Moreover, notice that the Teichmüller space defined in \cite{Nielsen} is connected but we use the definition of Teichmüller space \(\mathcal{T}_{Ein}\) given in \cite{looijenga2021teichmuller}, according to which it has multiple components  in general.

\begin{lem}\label{free_action}
	If \(X\) is such that \(\Aut(X)\rightarrow \bO^+(\bLambda_X)\) is injective then, for any component \(\mathcal{C}\) of \(\mathcal{T}_{Ein}\), we have \(\Mod(X)_{\mathcal{C}}\cap\T(X)=\{\id\}\). In particular, \(\T(X)\) acts freely on \(\mathcal{T}_{Ein}\) and hence the projection \(\mathcal{T}_{Ein}\rightarrow\mathcal{T}_{Ein}/\T(X)=:\mathcal{M}_{Ein}\) is a principal \(\T(X)\)-bundle.
\end{lem}
\begin{proof}
	Suppose \([\varphi]\in\Mod(X)_{\mathcal{C}}\cap\T(X)\). Since \([\varphi]\) acts trivially in cohomology, it fixes the period of any \([g]\in\mathcal{C}\) and it acts as the identity on the \(2\)-sphere of complex structures for \(g\). As a consequence of the equivariance of the period map and the fact that by \cite[Corollary 5.4]{looijenga2021teichmuller} the restriction of the period map to \(\mathcal{C}\) is an isomorphism onto its image, we have that \([\varphi]\) fixes the component \(\mathcal{C}\) point-wise. Hence, there is a representative \(\varphi\) that preserves an Einstein metric and a complex structure, and then \(\varphi\in\Aut(X)\). By hypothesis, the fact that \(\varphi\) acts trivially in cohomology implies that \(\varphi=\id\).
	
	 The second assertion follows by the fact that a non-trivial element of \(\T(X)\) does not fix any component, together with \(\T(X)\) being discrete. 
\end{proof}
We set
\(
\M := \mathcal{T}_{Ein} \times_{\Mod(X)} \E\Mod(X)\).
\begin{lem}\label{homotopy_equiv}
	Suppose \(X\) is such that \(\Aut(X)\rightarrow \bO^+(\bLambda_X)\) is injective,
	then there is a homotopy equivalence 
	\[
	\M \cong \mathcal{M}_{Ein} \times_{\Gamma} \E\Gamma.
	\]
\end{lem}
\begin{proof}
	Since \(\E\Mod(X)\times\E\Gamma\) is contractible and \(\Mod(X)\) acts freely on it, we have 
	\[\M\cong \mathcal{T}_{Ein}\times_{\Mod(X)}(\E\Mod(X)\times\E\Gamma),\] which has a fibration over \(\mathcal{T}_{Ein}\times_{\Mod(X)}\E\Gamma\) with contractible fiber \(\E\Mod(X)\). The fibration is a homotopy equivalence and since the action of \(\T(X)\) is free on \(\mathcal{T}_{Ein}\) by \autoref{free_action}, then \(\mathcal{T}_{Ein}\times_{\Mod(X)}\E\Gamma= \mathcal{M}_{Ein}\times_{\Gamma}\E\Gamma\).
\end{proof}
We are now ready to prove the first generalization.
\begin{proof}[Proof of \autoref{thm1}]
	Firstly, observe that by the description of the projection \(\mathcal{T}_{Ein}\rightarrow\mathcal{M}_{Ein}\) in \autoref{free_action}, any component of \(\mathcal{M}_{Ein}\) is diffeomorphic to any component of \(\mathcal{T}_{Ein}\) and hence \(\pi_1(\mathcal{T}_{Ein})=\pi_1(\mathcal{M}_{Ein})=1\).
	
	The long exact sequence of homotopy groups associated to the fibration $\mathcal{M}_{Ein} \to \M \to \B\Gamma$ implies that
	\[
	1 = \pi_1(\mathcal{M}_{Ein}) \to \pi_1(\M) \to \pi_1(\B\Gamma) \to \pi_0(\mathcal{M}_{Ein}) 
	\]
	is exact.
	The natural projection map $\M \to\B\Gamma$ induces an injection $\pi_1(\M) \hookrightarrow \pi_1(\B\Gamma)= \Gamma$, whose image consists of elements which are sent to the same point via the map \(\pi_1(\B\Gamma) \to \pi_0(\mathcal{M}_{Ein})\). The image of the latter is given by the orbit of a component \(\widetilde{\mathcal{C}}\) of \(\mathcal{M}_{Ein}\) for the monodromy action of \(\pi_1(\B\Gamma)= \Gamma\), hence \(\pi_1(\M)\) is isomorphic to the image via \(\rho\colon \Mod(X)\rightarrow \Gamma\) of the stabilizer of the component \(\widetilde{\mathcal{C}}\). If \(\mathcal{C}\) is a component of \(\mathcal{T}_{Ein}\) which is sent to \(\widetilde{\mathcal{C}}\), then the stabilizer of \(\widetilde{\mathcal{C}}\) corresponds to \(\Mod(X)_{\mathcal{C}}/(\Mod(X)_{\mathcal{C}}\cap \T(X))\) which equals \(\Mod(X)_{\mathcal{C}}\) by virtue of \autoref{free_action}. This shows that \(\pi_1(\M)\cong\Mon^2(X)\). 
	
	From the description of \(\M\), the map $\M \to \B\Gamma$ must factor as
	\[
	\M \to \B\Mod(X) \to \B\Gamma.
	\]
	Hence, the induced map $s \colon\pi_1(\M) \to \pi_1(\B\Mod(X)) = \pi_0(\Mod(X)) = \Mod(X)$ is a splitting of the restriction of $\rho\colon \Mod(X) \to \Gamma\subseteq \bO^+(\bLambda_X)$ to the stabilizer \(\Mod(X)_{\mathcal{C}}\).
\end{proof}

\begin{proof}[Proof of {\autoref{cor:max_monodromy}}]
	The map \(\Aut(X)\rightarrow \bO^+(\bLambda_X)\) is known to be injective by \cite[Lemma 3]{beauville1983varietes} for \(X\) of K3\(^{[n]}\) type, and by \cite[Theorem 2.1]{mongardi2014automorphisms} for \(X\) of OG10 type. The monodromy is maximal, i.e. \(\Mon^2(X)=\Gamma=\bO^+(\bLambda_X)\), by \cite[Lemma 9.2]{markman2011survey} for \(X\) of K3\(^{[n]}\) type with \(n-1\) a prime power and by \cite[Theorem 5.4]{onorati2020monodromy} for \(X\) of OG10 type.
\end{proof}

\begin{question}
	What can be said about the other known deformation types? 
\end{question}
The group \(\Mon^2(X)\) is available for all the known deformation types. If it is maximal, then \(\Gamma=\bO^+(X)\), but if it is a proper subgroup of \(\bO^+(\bLambda_X)\), then \(\Gamma\) is not known by the author.\\
Moreover, in the case \(\Aut(X)\rightarrow \bO^+(\bLambda_X)\) is not injective the argument given does not work, but it is a priori not clear if the same result might hold or not.

\section{Nielsen realization}

Let $G$ be a finite subgroup of $\Mod(X)$, by abuse of notation its image in $\bO^+(\bLambda_X)$ will be sometimes denoted again by $G$.
We want to give an answer to the Nielsen realization problem in terms of a \(G\)-invariant lattice, as done in \cite[Theorem 2.1]{farb2021nielsen} for K3 surfaces.
\\

Since $\Gr^+(3,\bLambda_X\otimes\mathbb{R})$ is the symmetric space of $\bO^+(\bLambda_X\otimes\mathbb{R})$, it is non-positively curved then $G$ must fix a point $P$. This means that $P$ is a $G$-invariant positive $3$-plane, and hence there is a linear representation $G\rightarrow \SO(P)$ of $G$.\\
Let $\II_G$ be the sum of all the irreducible $G$-subrepresentations of $\bLambda_X\otimes\mathbb{R}$ which are isomorphic to any of the ones appearing in $P$.  
\begin{definition}\label{def:LG}
	Let $ \bL_G=\II_G^\perp\cap \bLambda_X$.
\end{definition}
Notice that the same notation is sometimes used for the coinvariant lattice, but in fact the two lattices differ in general, for example \(\bL_G\) is always negative-definite but the coinvariant lattice of a non-symplectic automorphism is not definite.
\\

We say that a linear form $\delta\in\bLambda_X^\vee$ is \textit{negative} if its kernel has signature $(3,b_2(X)-4)$, or equivalently if its image via the embedding \(\bLambda_X^\vee\subset\bLambda_X\otimes\mathbb{Q} \) has negative square.
If $\mathcal{C}$ is a connected component of the Teichm\"uller space, let $\Delta_\mathcal{C}\subset \bLambda_X^\vee$ be the set of indivisible negative forms which are represented by an irreducible rational curve for a \hk{} metric belonging to $\mathcal{C}$. 

\begin{proof}[Proof of \autoref{thm3}]
	From the description in \cite[Section 5]{looijenga2021teichmuller}, each connected component $\mathcal{C}$ of the Teichm\"uller space is mapped diffeomorphically onto $$\Gr^+(3,\bLambda_X\otimes\mathbb{R})_{\Delta_\mathcal{C}}=\Gr^+(3,\bLambda_X\otimes\mathbb{R})- \bigcup_{\delta\in\Delta_\mathcal{C}}\Gr^+(3, \delta^\perp\otimes\mathbb{R})$$
	which is connected (and simply connected).
	This in particular means that if $G$ comes from a group of isometries for an Einstein metric, then the image $P$ via the period map is $G$-invariant and not orthogonal to any $\delta$, hence $\bL_G$ does not contain any $\delta$. If $G$ preserves a metric which is also Kähler, then the positive cone must be preserved by $G$ and we can find a $G$-invariant Kähler class which spans the trivial representation in $\bL_G^\perp$. For a detailed description of the period map we refer to \cite[Section 2.1]{looijenga2021teichmuller}.\\
	Suppose now that $G$ is a subgroup of $\Mod(X)$ which preserves a connected component of the Teichm\"uller space and for which $\bL_G$ does not contain any element in $\Delta_\mathcal{C}$. We argue as in the proof of \cite[Theorem 1.2]{farb2021nielsen}: among the $G$-invariant $3$-planes $P\subset \bL_G^\perp\otimes\mathbb{R}$, the ones such that $P^\perp\cap\bLambda_X=\bL_G$ are dense, so we can find a positive-definite $P\subset \bL_G^\perp\otimes\mathbb{R}$ such that $P^\perp\cap\bLambda_X=\bL_G$. Now, since $P$ does not lie in any $\delta^\perp$ for $\delta\in\Delta_\mathcal{C}$, the description of the image of the period map ensures that $G$ lifts to a group of isometries for a metric $g$ provided that it preserves a component $\mathcal{C}$ of the Teichm\"uller space. 
	Lastly, having the trivial representation in $\bL_G^\perp$ means that $G$ fixes a positive class $0\neq k\in P$ and hence the orientation determines a complex structure on $k^\perp \subset P$ which, again by surjectivity of the period map, is achieved by a complex structure on $X$ that makes $g$ a Kähler metric. \\
	The trivial representation is spanned by a positive integral $(1,1)$-class, so we can conclude using Huybrechts' projectivity criterion.

\end{proof}

The situation could be much more complicated than for K3 surfaces: as already noticed in \cite[Question 10.5]{markman2011survey} the stabilizer $\Mod(X)_\mathcal{C}$ could depend on the component $\mathcal{C}$ and it could intersect non-trivially the Torelli group, so it could happen that not every subgroup of $\bO^+(\bLambda_X)$ is the image of some stabilizer of a component and even those which are could have elements acting trivially on $\bLambda_X$. In case $G\subseteq \Mon^2(X)$ is the image of a group which intersects non-trivially $\T(X)$, then a lift could be found but it would be a finite extension of $G$.

\bibliographystyle{abbrv}

\bibliography{references}
\end{document}